\documentclass[11pt,a4paper]{article}
\usepackage{epsf,epsfig,amsfonts,amsgen,amsmath,amssymb,amstext,amsbsy,amsopn,amsthm,lineno}
\usepackage{color}
\usepackage{graphicx,tikz,caption,subfig,epstopdf}
\usepackage{enumerate}
\newtheorem{theorem}{Theorem}

\newtheorem{cor}{Corollary}

\theoremstyle{definition}

\newtheorem{case}{Case}

\definecolor{DarkGreen}{rgb}{0.2, 0.6, 0.3}

\tikzstyle{vertex}=[circle, draw, inner sep=0pt, minimum size=3pt]
\newcommand{\vertex}{\node[vertex]}

\def\qed{\hfill$\Box$\vspace{12pt}}

\begin{document}
\title
{\bf\ On characterizing the critical graphs for matching Ramsey numbers\thanks{The first author is supported by NSFC (No. 11701441). The third author is supported by NSFC (Nos. 11671320 and U1803263) and the Fundamental Research Funds for the Central Universities (No. 3102019ghjd003).}
}

\date{}

\author{\quad Chuandong Xu $^{a,}$\thanks{Corresponding author.},
\quad Hongna Yang $^{b,c}$,
\quad Shenggui Zhang $^{b,c,}$\thanks{E-mail addresses: xuchuandong@xidian.edu.cn (C.~Xu), yanghn@mail.nwpu.edu.cn (H.~Yang),
	sgzhang@nwpu.edu.cn (S.~Zhang).
}\\[2mm]
\small $^a$School of Mathematics and Statistics, Xidian University,\\
\small Xi'an, Shaanxi 710126, China\\
\small $^{b}$ School of Mathematics and Statistics, \\
\small Northwestern Polytechnical University, Xi'an, Shaanxi 710129, China\\
\small $^{c}$ Xi'an-Budapest Joint Research Center for Combinatorics, \\
\small Northwestern Polytechnical University, Xi'an, Shaanxi 710129, China
}
%\date{\today}
\maketitle

\begin{abstract}
Given simple graphs $H_{1},H_{2},\ldots,H_{c}$, the Ramsey number $r(H_{1},H_{2},\ldots,H_{c})$ is the smallest positive integer $n$ such that every edge-colored $K_{n}$ with $c$ colors contains a subgraph in color $i$ isomorphic to $H_{i}$ for some $i\in\{1,2,\ldots,c\}$. The critical graphs for $r(H_1,H_2,\ldots,H_c)$ are edge-colored complete graphs on $r(H_1,H_2,\ldots,H_c)-1$ vertices with $c$ colors which contain no subgraphs in color $i$ isomorphic to $H_{i}$ for any $i\in \{1,2,\ldots,c\}$. For $n_1\geq n_2\geq \ldots\geq n_c\geq 1$, Cockayne and Lorimer (The Ramsey number for stripes, {\it J.\ Austral.\ Math.\ Soc.} \textbf{19} (1975), 252--256.) showed that $r(n_{1}K_{2},n_{2}K_{2},\ldots,n_{c}K_{2})=n_{1}+1+ \sum\limits_{i=1}^c(n_{i}-1)$, in which $n_{i}K_{2}$ is a matching of size $n_{i}$. Using the Gallai-Edmonds Theorem, we characterized all the critical graphs for $r(n_{1}K_{2},n_{2}K_{2},\ldots,n_{c}K_{2})$, implying a new proof for this Ramsey number.
\medskip

\smallskip

\end{abstract}

\section{Introduction}
All graphs considered in this paper are finite and simple. For terminology and notation not defined here, we refer the reader to Bondy and Murty \cite{Bondy}.

An edge-colored graph is {\it monochromatic} if all its edges have the same color. Given simple graphs $H_{1},H_{2},\ldots,H_{c}$, the {\it Ramsey number} $r(H_1,H_2,\ldots,H_c)$ is the  smallest positive integer $n$ such that every {\it $c$-edge-coloring} of $K_{n}$ (an assignment of $c$ colors to the edges of $K_{n}$) contains a monochromatic subgraph in some color $i\in\{1,2,\ldots,c\}$ isomorphic to $H_i$. A {\it critical graph} for $r(H_{1},H_{2},\ldots,H_{c})$ is a $c$-edge-colored complete graph on $r(H_{1},H_{2},\ldots,H_{c})-1$ vertices, which contains no subgraphs in color $i$ isomorphic to $H_{i}$ for any $i\in \{1,2,\ldots,c\}$. 

Determining the value of classical Ramsey numbers seems to be extremely hard (see \cite{Rad2017} for a survey). But for multiple copies of graphs, Burr, Erd{\H o}s and Spencer \cite{BES1975} obtained surprisingly sharp and general upper and lower bounds on $r(nG,nH)$ for fixed $G$, $H$ and sufficiently large $n$. They also showed that $r(mK_3,nK_3)=3m+2n$ when $m\geq n$, $m\geq 2$. Hook and Isaak \cite{HI} made a conjecture on the critical graphs for $r(mK_3,nK_3)$. Another well-known result in this area is due to Cockayne and Lorimer \cite{CL}. 
\begin{theorem}[Cockayne and Lorimer \cite{CL}]
For $n_1\geq n_2\geq \ldots\geq n_c\geq 1$,
	\[
		r(n_{1}K_{2},n_{2}K_{2},\ldots,n_{c}K_{2})=n_{1}+1+ \sum\limits_{i=1}^c(n_{i}-1).
	\]
\end{theorem}
This result has been generalized to complete graphs versus matchings by Lorimer and Solomon \cite{LorSol1992}, and to hypergraphs by Alon et al.~\cite{AF}. For the Ramsey number of matchings, Hook and Isaak \cite{HI} characterized the critical graphs for $r(mK_2,nK_2)$ for $m\geq n\geq 1$. The class of all critical graphs for $r(n_{1}K_{2},n_{2}K_{2},\ldots,n_{c}K_{2})$ has not been determined yet.

Cockayne and Lorimer \cite{CL} gave a critical graph for $r(n_{1}K_{2},n_{2}K_{2},\allowbreak\ldots,n_{c}K_{2})$ which is a $c$-edge-colored complete graph $G$ on $n_{1}+ \sum\limits_{i=1}^c(n_{i}-1)$ vertices whose vertex set $V(G)$ has $c$ parts $V_{1},\ldots,V_{c}$ such that $\left|V_{1}\right|=2n_{1}-1$, $\left|V_{i}\right|=n_{i}-1$ for $i\in \{2,\ldots,c\},$ and the color of an edge $e=xy$ in $G$ is the maximum $j$ for which $\{x, y\}$ has a non-empty intersection with $V_{j}$. It is easy to see that $G$ contains no monochromatic $n_{i}K_{2}$ in color $i$ for any $i\in \{1,2,\ldots,c\}$.

\begin{figure}[bht]
\[\begin{tikzpicture}
\draw (0,0.5) ellipse   (0.7 and 0.5);
\draw (-3,-2) ellipse   (0.6 and 0.4);
\draw (-0.5,-2) ellipse (0.6 and 0.4);
\draw (3,-2) ellipse    (0.6 and 0.4);

\node (n1) at (0,0.5)      {\small{$K_{2n_{1}-1}$}};
\node (n2) at (-3,-2)    {\footnotesize{$K_{n_{2}-1}$}};
\node (n3) at (-0.5,-2)     {\footnotesize{$K_{n_{3}-1}$}};
\node (n4) at (3,-2)     {\footnotesize{$K_{n_{c}-1}$}};
\node (n5) at (1.25,-1.95) {$\cdots$ };

\draw[semithick](-0.67,0.38)--(-2.9,-1.6)
node [midway, above] {\small{2}};
\draw [semithick](-0.5,0.15)--(-2.6,-1.7)
node [midway, below] {\small{2}};

\draw [thick,dash pattern=on 7pt off 4pt](0,0)--(-0.25,-1.65)
node [midway, right=-3] {\small{3}};
\draw [thick,dash pattern=on 7pt off 4pt](-0.27,0.05)--(-0.53,-1.6)
node [midway, left=-3] {\small{3}};

\draw [thick,dash pattern=on 1pt off 2.5pt on 1pt off 2.5pt](0.5,0.1)--(n4)
node [midway, below] {$c$};
\draw [thick,dash pattern=on 1pt off 2.5pt on 1pt off 2.5pt](0.68,0.3)--(2.96,-1.61)
node [midway, above] {$c$};

\draw [semithick](-2.4,-1.9)--(-1.1,-1.9)
node [midway, above=-2.7] {\small{2}};
\draw [thick,dash pattern=on 7pt off 4pt](-2.4,-2.1)--(-1.1,-2.1)
node [midway, below=-2.5] {\small{3}};

\draw (0.08,-2)edge[style={bend right}, thick,dash pattern=on 7pt off 4pt](2.4,-2)
node (n6) at (1.25,-2.18) {\small{3}};
\draw (0.08,-2.2)edge[style={bend right}, thick,dash pattern=on 1pt off 2.5pt on 1pt off 2.5pt](2.5,-2.2)
node (n6) at (1.25,-2.7) {$c$};

\draw (-2.7,-2.35)edge[style={bend right}, semithick](2.7,-2.35)
node (n6) at (0,-2.98) {\small{2}};
\draw (-3.15,-2.4)edge[style={bend right}, thick,dash pattern=on 1pt off 2.5pt on 1pt off 2.5pt](3.15,-2.4)
node (n6) at (0,-3.46) {$c$};

\end{tikzpicture}\]
\caption{The structure of the critical graphs for $r(n_1K_1,n_2K_2,\ldots,n_cK_2)$.}
\label{fig:math:G}
\end{figure}
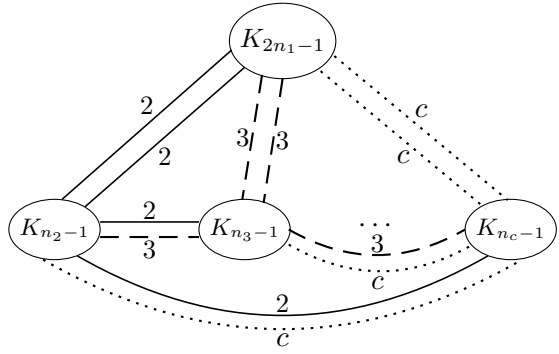

Motivated by Cockayne and Lorimer's result, in this paper we studied the structure of the critical graphs for $r(n_{1}K_{2},n_{2}K_{2},\ldots,n_{c}K_{2})$ (see Figure \ref{fig:math:G} for an example).

\begin{theorem}
\label{th 1}
For $n_{1}\geq n_{2}\geq \ldots\geq n_{c}\geq 1$, let $G$ be a $c$-edge-colored complete graph with order $n\geq n_{1}+ \sum\limits_{i=1}^c(n_{i}-1)$. If $G$ contains no monochromatic $n_{i}K_{2}$ in color $i$ for any $i\in \{1,2,\ldots,c\}$, then $n=n_{1}+ \sum\limits_{i=1}^c(n_{i}-1)$ and the colors of $G$ can be relabeled such that:
\begin{enumerate}
 \item [(a)] $V(G)$ can be partitioned into $c$ parts $V_{1},V_2, \ldots,V_{c}$, where $\left|V_{1}\right|=2n_{1}-1$, $\left|V_{i}\right|=n_{i}-1$, and all the edges with ends both in $V_{i}$ have color $i$, for $i\in\{1, 2, \ldots, c\}$;
 \item [(b)] all the edges with one end in $V_{1}$ and the other end in $V_{i}$ have color $i$, for $i\in\{2,\ldots,c\}$;
 \item [(c)] all the edges with one end in  $V_{i}$ and the other end in $V_{j}$ have color either $i$ or $j$, for $\{i,j\}\subseteq\{2,\ldots,c\}$.
\end{enumerate}
\end{theorem}

Bialostocki and Gy\'{a}rf\'{a}s \cite{BG} showed that Cockayne and Lorimer's proof (there is a gap, a missed case, in this proof) can be modified to give a more general result.

\begin{theorem}[Bialostocki and Gy\'{a}rf\'{a}s \cite{BG}]
	for $n_1\geq n_2\geq \ldots\geq n_c\geq 1$ and $n\geq n_1+1+\sum\limits_{i=1}^c(n_{i}-1)$, every $c$-edge-colored $n$-chromatic graph contains a monochromatic $n_{i}K_{2}$ for some $i\in\{1,2,\ldots c\}$.
\end{theorem}

As mentioned in \cite{BG}, Zolt\'an Kir\'aly pointed out that the $n$-chromatic graph version result can be deduced from the complete graph version result. Here we will show that Zolt\'an Kir\'aly's method can work for more general graph classes. Let $G$ be an edge-colored graph with $c$ colors. If there is a partition $\{V_{1},V_{2},\ldots,V_{n}\}$ of $V(G)$ such that $E(V_{i}, V_{j})\neq \emptyset$ for $i\neq j$ and $n\geq n_{1}+1+\sum\limits_{i=1}^c(n_{i}-1)$, then by identifying each $V_i$ to a single vertex $v_i$ and deleting the multiplied edges, one can obtain a $c$-edge-colored complete graph on $n$ vertices, denoted by $G^*$. It's easy to see that each monochromatic $n_iK_2$ with some color $i$ in $G^*$ corresponds to a monochromatic $n_iK_2$ with color $i$ in $G$. 
\begin{cor}
\label{ob 1}
Let $G$ be an edge-colored graph with $c$ colors. If there is a partition $\{V_{1},V_{2},\ldots,V_{n}\}$ of $V(G)$ such that $E(V_{i}, V_{j})\neq \emptyset$ for each $i\neq j$ and $n\geq n_{1}+1+\sum\limits_{i=1}^c(n_{i}-1)$, then $G$ contains a monochromatic $n_{i}K_{2}$ for some $i\in\{1,2,\ldots,c\}$.
\end{cor}

The proof of Theorem \ref{th 1} is in Section 2. At the end of this paper, we remark a simple application of Theorem \ref{th 1}.

\section{Proof of Theorem 1}

First, we will state the Gallai-Edmonds Theorem which plays an essential role in our proof.

Let $M$ be a matching of a graph $G$ with order $n$. Each vertex incident with an edge in $M$ is said to be {\it covered} by $M$. A {\it maximum matching} of $G$ is a matching that covers as many vertices as possible. When $n$ is even (odd), a {\it perfect matching} ({\it near-perfect matching}) is a maximum matching of $G$ which covers $n$ vertices ($n-1$ vertices). We call $G$ {\it factor-critical} if $G-v$ has a perfect matching for each vertex $v\in G$.

For a graph $G$, let $D(G)$ be the set of vertices that cannot be covered by at least one maximum matching of $G$, $A(G)$ be the set of vertices that have neighbours in $D(G)$, and $C(G)=V(G)\setminus (D(G)\cup A(G))$. The following Gallai-Edmonds Theorem is due to Gallai \cite{G2} and Edmonds \cite{E1}. The current version of this theorem we used here can be found in Lov{\'a}sz and Plummer \cite{LP} (pp. 94, Theorem 3.2.1). We call $D(G)$, $A(G)$, and $C(G)$ the {\it Gallai-Edmonds decomposition} of $G$ (see Figure \ref{fig:match:CAD} as an example).
\begin{theorem}[Gallai-Edmonds Theorem]
\label{GallaiEdmonds}
 For a graph $G$, let $D(G)$, $A(G)$, and $C(G)$ be defined as above. Then
\begin{enumerate}
\item [(a)] the components of the subgraph induced by $D(G)$ are factor-critical;
\item [(b)] the subgraph induced by $C(G)$ has a perfect matching;
\item [(c)] the bipartite graph obtained from $G$ by deleting the vertices of $C(G)$ and the edges spanned by $A(G)$ and by contracting each component of $D(G)$ to a single vertex has a positive surplus (as viewed from $A(G)$, i.e., $|N(S)| - |S| > 0$ for each nonempty subset $S$ of $A(G)$);
\item [(d)] if $M$ is any maximum matching of $G$, it contains a near-perfect matching of each component of $D(G)$, a perfect matching of each component of $C(G)$ and matches all vertices of $A(G)$ with vertices in distinct components of $D(G)$;
\item [(e)] the size of a maximum matching $M$ is equal to $\frac{1}{2}(\left| V(G)\right|-\omega(D(G))+ \left| A(G)\right|)$, where $\omega(D(G))$ denotes the number of components of the graph spanned by $D(G)$.
\end{enumerate}
\end{theorem}

Since there exists no monochromatic $n_{i}K_{2}$ in color $i$ in color class $G^i$ (the subgraph of $G$ induced by all the edges in color $i$) for each $i\in \{1,2,\ldots,c\}$, we know that the {\it matching number} (the size a maximum matching) of $G^i$ is at most $n_{i}-1$. The Gallai-Edmonds Theorem characterizes the structure of a graph based on its matching number. We will deduce from the Gallai-Edmonds Theorem that each color class $G^i$ in $G$ cannot have too many edges. On the other hand, the union of these color classes have to cover all the edges of $G$. Finally we characterize the structure of $G$, which also implies a new proof on the value of  $r(n_{1}K_{2},n_{2}K_{2},\ldots,n_{c}K_{2})$.

\vskip \baselineskip

\begin{Tproof}\textbf{ of Theorem \ref{th 1}.}
Suppose that $G$ has $n\geq n_{1}+ \sum\limits_{i=1}^c(n_{i}-1)$ vertices and contains no monochromatic $n_{i}K_{2}$ in color $i$ for any $i\in \{1,2,\ldots,c\}$. If $n_i=1$ for some $1\leq i\leq c$, then $G$ contains no edges with color $i$. We can ignore color $i$ in our discussion and there is no influence to the conclusions. So we will assume $n_{1}\geq n_{2}\geq \ldots\geq n_{c}\geq 2$ in this proof.

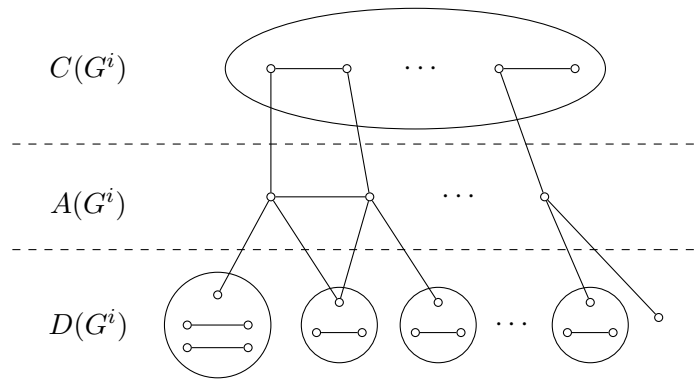
\begin{figure}[bht]
	\[\begin{tikzpicture}
	%~~~~~~~~~~~~~~~~~~~~~~~~~~~draw grids ~~~~~~~~~~~~~~~~~~
	%    \draw[step=1cm,gray,very thin] (-5,2) grid (5,4);
	%~~~~~~~~~~~~~~~~~~~~~~~~~~~~~~~~~~~~~~~~~~~~~~~
	\draw (-2.8,-0.6) ellipse (0.7 and 0.7) ;
	\draw (-1.2,-0.6) ellipse (0.5 and 0.5) ;
	\draw (0.1,-0.6) ellipse (0.5 and 0.5) ;
	\draw (2.1,-0.6) ellipse (0.5 and 0.5) ;
	\draw (-0.2,2.8) ellipse (2.5 and 0.8) ;

	\vertex (u1) at (-2.1,2.8) {};
	\vertex (u2) at (-1.1,2.8) {};
	\vertex (u3) at (0.9,2.8) {};
	\vertex (u4) at (1.9,2.8) {};
	
	\vertex (v1) at (-2.1,1.1) {};
	\vertex (v2) at (-0.8,1.1) {};
	\vertex (v3) at (1.5,1.1) {};
	
	\vertex (w1) at (-2.8,-0.2) {};
	\vertex (w2) at (-3.2,-0.6) {};
	\vertex (w3) at (-2.4,-0.6) {};
	\vertex (w4) at (-3.2,-0.9) {};
	\vertex (w5) at (-2.4,-0.9) {};
	\vertex (w6) at (-1.5,-0.7) {};
	\vertex (w7) at (-0.9,-0.7) {};
	\vertex (w8) at (-1.2,-0.3) {};
	\vertex (w9) at (2.1,-0.3) {};
	\vertex (w11) at (2.4,-0.7) {};
	\vertex (w12) at (1.8,-0.7) {};
	\vertex (w10) at (3,-0.5) {};
	\vertex (x1) at (-0.2,-0.7) {};
	\vertex (x2) at (0.4,-0.7) {};
	\vertex (x3) at (0.1,-0.3) {};
	
	\node at (0.4,1.1) {$\cdots$};
	\node at (-0.1,2.8) {$\cdots$};
	\node at (1.1,-0.6) {$\cdots$};
	\node at (-4.5,2.8) {$C(G^i)$};
	\node at (-4.5,1) {$A(G^i)$};
	\node at (-4.5,-0.6) {$D(G^i)$};
	
	\path (u1) edge (u2) (u3) edge (u4) (u1) edge (v1) (u2) edge (v2) (u3) edge (v3) (w12) edge (w11) (v1) edge (v2)
	(v1) edge (w1) (v1) edge (w8)(v2) edge (w8) (v3) edge (w9)(v2) edge (x3) (x1) edge (x2)(v3) edge (w10)
	(w2) edge (w3) (w4) edge (w5) (w6) edge (w7);

	\draw[dashed] (-5.5,1.8) -- (3.5,1.8);
	\draw[dashed] (-5.5,0.4) -- (3.5,0.4);
	
	\end{tikzpicture}\]
	\caption{The Gallai-Edmonds decomposition of the color class $G^i$.}
	\label{fig:match:CAD}
\end{figure}

Let $G^1, G^2, \ldots, G^c$  be the color classes of $G$. For each $i\in\{1,2,\ldots,c\}$, the matching number of $G^i$ is at most $n_{i}-1$ since $G$ contains no monochromatic $n_{i}K_{2}$ in color $i$. Let $C(G^i)$, $A(G^i)$, and $D(G^i)$ be the Gallai-Edmonds decomposition of $G^i$ (see Figure \ref{fig:match:CAD}). Denote the vertex sets of components in $G^i[D(G^i)]$ by $D_1(G^i), D_2(G^i) \ldots D_{t_i}(G^i)$. Let 
\[
	a_{i}=\left| A(G^i)\right|,\quad d_{i_{0}}=\frac{\left| C(G^i)\right|}{2},\quad\ d_{i_{k}}=\frac{\left| D_{k}(G^i)-1\right|}{2}\ \text{for}\ k\in\{1,2,\ldots,t_i\}.
\]
By the Gallai-Edmonds Theorem,  $a_{i}+d_{i_0}+d_{i_1}+\cdots+d_{i_{t_i}}$ is the mathcing number of $G^i$. Since the matching number of $G^i$ is at most $n_i-1$, there holds

\begin{equation*}
d_{i_{0}}+d_{i_{1}}+\ldots+d_{i_{t_i}}\leq n_{i}-1-a_{i}.
\end{equation*}

The following inequalities give an upper bound on the number of edges with its ends both in $C(G^i)$ or in $D(G^i)$, in which the third inequality can be checked by comparing the size of a complete graph with order $2(d_{i_{0}}+d_{i_{1}}+\cdots+d_{i_{t_i}})+1$ and the size of a subgraph of it. We have
\begin{equation}\label{ineq_1}
	\begin{aligned}
	\left|E(G^i[C(G^i)])\right|+\left|E(G^i[D(G^i)])\right|
	&\leq\binom{2d_{i_{0}}}{2}+\binom{2d_{i_{1}}+1}{2}+\cdots+\binom{2d_{i_{t_i}}+1}{2}\\
	&\leq \binom{2d_{i_{0}}+1}{2}+\binom{2d_{i_{1}}+1}{2}+\cdots+\binom{2d_{i_{t_i}}+1}{2}\\
	&\leq \binom{2(d_{i_{0}}+d_{i_{1}}+\cdots+d_{i_{t_i}})+1}{2}\\
	&\leq \binom{2(n_{i}-1-a_{i})+1}{2}.
	\end{aligned}
\end{equation}

Next, we give bounds on the number of edges incident with vertices in $A(G^i)$ which can be partitioned into $a_{i}$ stars. There are $\sum \limits_{i=1}^c a_{i}$ such stars in total. Let $H$ be the subgraph of $G$ with vertex set $V(G)$ and edge set the union of the edge sets of these stars. Those vertices in $V(G)-\cup_{i=1}^{c}A(G^i)$ form an independent set of size at least $ n-\sum \limits_{i=1}^c a_{i}$ in $H$. Thus $H$ has at most $\binom{n}{2}- \binom{n-\sum a_{i}}{2}$ edges. Together with the edges in $G^i[C(G^i)]$ and $G^i[D(G^i)]$ for $1\leq i\leq c$, we have an upper bound on the number of edges in $\cup_{i=1}^{c}G^i$ which is a complete graph with oder $n$:
\begin{equation}\label{ineq_2}
\binom{n}{2}-\binom{n-\sum \limits_{i=1}^c a_{i}}{2}+\sum \limits_{i=1}^c\binom{2(n_{i}-1-a_{i})+1}{2}\geq \binom{n}{2}.
\end{equation}

Note that $n\geq n_{1}+\sum \limits_{i=1}^c(n_{i}-1-a_{i})$. There follows
\begin{equation}\label{ineq_3}
\sum \limits_{i=1}^c\binom{2(n_{i}-1-a_{i})+1}{2}\geq \binom{n-\sum \limits_{i=1}^c a_{i}}{2}\geq \binom{n_{1}+\sum \limits_{i=1}^c(n_{i}-1-a_{i})}{2}.
\end{equation}
For the convenience of discussion, let $b_{i}=n_{i}-1-a_{i}$ for $1\leq i\leq c$. Then we have
\begin{equation}\label{ineq_4}
\sum \limits_{i=1}^c\binom{2b_{i}+1}{2}\geq \binom{n_{1}+\sum \limits_{i=1}^c b_{i}}{2}.
\end{equation}

We will deduce the structure of $G$ from the above inequality. Assuming $b_{m}=\max\{ b_{1}, b_{2}, \ldots, b_{c} \}$, we get $b_{m}>0$ (otherwise (\ref{ineq_4}) dosen't hold since $n_1\geq 2$) and $0\leq b_{i}\leq b_{m}\leq n_{1}-1$. For $b_{i}>0$ and $i\neq m$, there holds $b_i\leq n_1-1\leq n_1-1+n_1-2$, i.e., $\frac{b_i+3}{2}\leq n_1$. There holds
\begin{equation}\label{ineq_5}
\begin{aligned}
\binom{2b_{i}+1}{2}
&=\binom{b_{i}}{2}+b_{i}(b_{i}+1)+\binom{b_{i}+1}{2}\\
&= \binom{b_{i}}{2}+b_{i}\cdot b_{i}+b_{i}\cdot \frac{b_{i}+3}{2}\\
&\leq \binom{b_{i}}{2}+b_{i}\cdot b_{m}+b_{i}\cdot n_{1}.
\end{aligned}
\end{equation}
The equality in (\ref{ineq_5}) holds if and only if $b_{i}=b_{m}$ and $\frac{b_{i}+3}{2}=n_{1}$, which only holds when $n_{1}=2$ and $b_{i}=b_{m}=1$. 

The last inequality in the following can be checked by treating each item as the size of a subgraph of a complete graph with order $n_{1}+\sum \limits_{i=1}^c b_{i}$. It follows from (\ref{ineq_5}) that
\begin{equation}\label{ineq_6}
\begin{aligned}
\sum \limits_{i=1}^c \binom{2b_{i}+1}{2}
&=\binom{2b_{m}+1}{2}+\sum \limits_{i=1,i\neq m}^c\binom{b_{i}+1}{2}\\
&\leq \binom{n_{1}+b_{m}}{2}+\sum \limits_{i=1,i\neq m}^c \left[ \binom{b_{i}}{2}+b_{i}\cdot b_{m}+b_{i}\cdot n_{1} \right]\\
&\leq \binom{n_{1}+\sum \limits_{i=1}^c b_{i}}{2}.
\end{aligned}
\end{equation}
The equalities in (\ref{ineq_6}) hold if and only if $b_{m}=n_{1}-1$ and there exists at most one nonzero $b_{i}$ with $i\neq m$. 

By (\ref{ineq_4}) and (\ref{ineq_6}), we get
$$\sum \limits_{i=1}^c\binom{2b_{i}+1}{2}= \binom{n_{1}+\sum \limits_{i=1}^c b_{i}}{2}.$$
Hence, the equalities hold throughout in inequalities (\ref{ineq_1})--(\ref{ineq_6}). Thus $n=n_{1}+\sum\limits_{i=1}^c(n_{i}-1)$ and $b_{m}=n_{1}-1$. Since $b_{m}=n_{m}-1-a_{m}$, $n_m\leq n_1$, and $a_m\geq 0$, there holds $n_m=n_1$ and $a_m=0$. Hence we can switch the colors of $G^1$ and $G^m$ to set $m=1$. There are two cases for the values of $b_{1}, b_{2}, \ldots,b_{c}$. 

\begin{case}
$b_{1}=n_{1}-1$, $b_{2}=\cdots=b_{c}=0.$
\end{case}

It follows that $a_{1}=0$, $a_{2}=n_{2}-1$, $\cdots$, $a_{c}=n_{c}-1$. For $i=1$, since the equality holds in inequality (1), there follows $C(G^1)=A(G^1)=\emptyset$ and $G^1[D(G^1)]\cong K_{2n_{1}-1}$. Thus $G^1\cong K_{2n_{1}-1}$. 

For $i\geq 2$, it follows from $a_{i}=n_{i}-1$ that $C(G^i)=\emptyset$, and components in $G^i[D(G^i)]$ are isolate vertices. Recall that $H$ contains the $a_i=n_i-1$ stars in color $i$, i.e., $H$ contains $G^i$.  Moreover,  $H=\bigcup\limits_{i=2}^{c}G^i\cong K_{n}\backslash E(K_{2n_{1}-1})$ (the complement of $K_{2n_1-1}$ in $K_n$). Thus $G$ has the required structure.
\begin{case}
$n_{1}=2$, $b_{1}=b_{2}=1$, $b_{3}=\cdots=b_{c}=0.$
\end{case}
It follows that $n_{1}=n_{2}=\cdots=n_{c}=2$ since $2\leq n_{i}\leq n_{1}$. For $i\neq 1$ and $b_{i}>0$, we assume $i=2$ for convenience. Thus $n=c+2$, $a_{1}=a_{2}=0$ and $a_{3}=\cdots=a_{c}=1$. By (1), $\left| E(G^1)\right| =\left| E(G^2)\right|=3$, and thus $G_{1}\cong G_{2}\cong K_{3}$. Also by (1), $G^1 \cup G^2$ is isomorphic to $K_{4}$, a contradiction.~
\qed
\end{Tproof}

\section{Remark}
Let $K_{n-1}\sqcup K_{1,k}$ be the graph obtained from $K_{n-1}$ by adding a new vertex $v$ and joining $v$ to $k$ vertices of $K_{n-1}$. For $n=r(H_{1}, H_{2},\ldots,H_{c})$, the {\it star-critical Ramsey number} is the smallest positive integer $k$ such that every $c$-edge-coloring of  $K_{n-1}\sqcup K_{1,k}$ contains a subgraph isomorphic to $H_{i}$ in color $i$ for some $i\in\{1,2,\ldots,c\}$, denoted by  $r_{*}(H_{1}, H_{2},\ldots,H_{c})$. This concept was introduced by Hook and Isaak \cite{HI}, who showed that $r_{*}(sK_{2}, tK_{2})=t$ for $s\geq t\geq 1$. The star-critical Ramsey numbers of other graphs have been investigated in \cite{HI,HLQ,HL,HM,LL,WS,ZB}.

A {\it $(H_{1}, H_{2},\ldots,H_{c})$-free coloring} of $K_{n-1}$ is a $c$-edge-coloring of $K_{n-1}$ that contains no subgraphs isomorphic to $H_{i}$ in color $i$ for any $i\in \{1,\ldots,c\}$. Thus every critical graph for $r(n_{1}K_{2},n_{2}K_{2},\ldots,n_{c}K_{2})$ has an  $(n_{1}K_{2},n_{2}K_{2},\ldots,n_{c}K_{2})$-free coloring.
By using Theorem 1, we get the following result on the star-critical Ramsey number of matchings.

\begin{theorem}
\label{the 0}
For $n_{1}\geq n_{2}\geq \ldots\geq n_{c}\geq 1$, let $r_{*}(n_{1}K_{2},n_{2}K_{2},\ldots,n_{c}K_{2})=\sum\limits_{i=2}^c(n_{i}-1)+1$.
\end{theorem}
\begin{Tproof}\textbf{.}
For convenience, we let $$n:=r(n_{1}K_{2},n_{2}K_{2},\ldots,n_{c}K_{2})=n_{1}+1+ \sum\limits_{i=1}^c(n_{i}-1),\quad m:=\sum\limits_{i=2}^c(n_{i}-1).$$ 

To show $r_{*}(n_{1}K_{2},n_{2}K_{2},\ldots,n_{c}K_{2})\geq m+1$, we give an $(n_{1}K_{2},n_{2}K_{2},\ldots,n_{c}K_{2})$-free coloring of $K_{n-1}\sqcup K_{1,m}$, which is constructed by a critical graph on $n-1$ vertices as defined in Theorem \ref{th 1} and a vertex $v$ with edges to each monochromatic $K_{n_{i}-1}$ colored by $i$ for $i\in \{2,\ldots,c\}$. 

Next we prove the reverse. Let $G$ be an edge-colored $K_{n-1}\sqcup K_{1,m+1}$ with $c$ colors, $H$ be the $K_{n-1}$ in $G$, and $v$ be the center of the star $K_{1,m+1}$. By Theorem \ref{th 1}, either $H$ contains a monochromatic $n_iK_2$ and we are done, or $H$ is a critical graph and contains an monochromatic $K_{2n_{1}-1}$ with some color, say color $1$. In the following we assume that $H$ belongs to the latter case. Thus no edges incident to $v$ has color $1$ in $G$, or there is a monochromatic $n_1K_2$. So the colors of the edges incident to $v$ belong to $\{2,\ldots,c\}$. Note that $n-1-m=2n_1$, there exists an edge $uv$ with $u\in V(H^1)$ ($H^1$ is the monochromatic $K_{2n_{1}-1}$ in H). Denote the color of $uv$ by $j$ ($j\in \{2,\ldots,c\}$). Then the edge $uv$ and an $(n_j-1)$-matching in $H^j$ form an $n_{j}K_{2}$ with color $j$ in $G$. The result follows.
\qed
\end{Tproof}

\end{document}